\documentclass[12pt]{article}
\usepackage{amssymb, amsmath}
\begin{document}
\voffset=0.0truein \hoffset=-0.5truein \setlength{\textwidth}{6.0in}
\setlength{\textheight}{8.8in} \setlength{\topmargin}{-0.2in}
\renewcommand{\theequation}{\arabic{section}.\arabic{equation}}
\newtheorem{thm}{Theorem}[section]
\newtheorem{lemma}{Lemma}[section]
\newtheorem{pro}{Proposition}[section]
\newtheorem{cor}{Corollary}[section]
\newcommand{\n}{\nonumber}
\newcommand{\w}{\omega}
\newcommand{\g}{\gamma}
\newcommand{\tv}{\tilde{v}}
\newcommand{\rb}{\bar{\rho}}
\newcommand{\tw}{\tilde{\omega}}
\renewcommand{\a}{\alpha}
\renewcommand{\o}{\omega}
\renewcommand{\O}{\Omega}
\newcommand{\e}{\varepsilon}
\renewcommand{\t}{\theta}
\newcommand{\nao}{\nabla ^\bot \theta}
\newcommand{\vare}{\varepsilon}
\newcommand{\bb}{\begin{equation}}
\newcommand{\ee}{\end{equation}}
\newcommand{\bq}{\begin{eqnarray}}
\newcommand{\eq}{\end{eqnarray}}
\newcommand{\bqn}{\begin{eqnarray*}}
\newcommand{\eqn}{\end{eqnarray*}}
\title{The finite time blow-up for the \\
Euler-Poisson equations in $\Bbb R^n$}
\author{Dongho Chae \\
Department of Mathematics\\
              Sungkyunkwan University\\
               Suwon 440-746, Korea\\
              {\it e-mail : chae@skku.edu}}
 \date{}
\maketitle
\begin{abstract}
We prove the finite time blow-up for $C^1$ solutions to the  Euler-Poisson equations in $\Bbb R^n$, $n\geq 1$, with/without background density for initial data satisfying suitable conditions. We also find a sufficient condition for the initial data such that $C^3$ solution  breaks down in finite time for the compressible Euler equations for polytropic gas flows.
\end{abstract}
\noindent{{\bf AMS subject classification:} 35Q35, 35B30}\\
\noindent{{\bf Key Words:} Euler-Poisson equations, finite time blow-up, polytropic flow}

\section{The Euler-Poisson equations}
 \setcounter{equation}{0}

We are concerned on the following (pressureless) Euler-Poisson equations in $\Bbb R^n$, $n\geq 1$.
\bb\label{ep1}
 \left\{ \aligned
 &\partial_t \rho+ \mathrm{div}\,(\rho v)=0,\\
 &\partial_t (\rho v)+\mathrm{div}\,(\rho u\otimes u) =-k\rho \nabla \phi,\\
&\Delta \phi =\rho,\\
& v(x,0)=v_0 (x), \,\, \rho(x,0)=\rho_0 (x).\\
\endaligned\right.
 \ee
The unknowns are the velocity field $v=(v^1,\cdots , v^n )=v(x,t),$ mass density $\rho=\rho(x,t) \geq 0$, and the potential function of the force field $\phi=\phi (x,t)$. We consider the case $k>0$, which  represents that the force is attractive.
Physically the system (\ref{ep1}) represents the large scale dynamics stars or cluster of stars, regarded as particles consisting of a gas, interacting through gravitation.
The global regularity/finite time blow-up problem from the view point of the `critical threshold phenomena' have
been studied extensively by Tadmor and his collaborators for the one dimensional Euler-Poisson equation  (with/without pressure term) (see e.g. \cite{eng1, eng2, tad} and references therein). In the multi-dimensional problem, the finite time blow-up for the attractive force case  with pressure
in the spherical symmetry  for compactly supported $\rho(x,t)$ in $\Bbb R^3$ was obtained in \cite{mak}, and for the repulsive case with similar geometry  the blow-up was deduced in \cite{per}.
On the other hand, the global regularity for some large class of initial data near a steady state is obtained in \cite{guo}. In \cite{smo} the stability type of result was obtained with inclusion of the pressure. As for the study for the equations (\ref{ep1}) in the general setting, which is our main concern  in this paper, the main difficulty was the forcing term $k\nabla \phi$, which is
of nonlocal nature,  and resembles the  notorious pressure term in the 3D incompressible Euler equations.
This feature was emphasized  in \cite{liu}, and was the main motivation for studying a model problem, called  `restricted Euler-Poisson equations', where the nonlocal forcing term is replaced by a local one. The finite time  blow-up for the restricted Euler equations was shown via spectral dynamics in \cite{liu}.
Our aim in this  paper is to show that the finite time blow-up is unavoidable for the original Euler-Poisson equations under suitable  conditions  for the initial data of large class.
At least for our blow-up problem the nonlocality can be bypassed simply by decomposing the evolution equations of the velocity gradient matrix into the symmetric part(deformation tensor equation) and the antisymmetric art(vorticity equation), and observing that the vanishing property of the vorticity is preserved along the particle trajectories, concentrating on the evolution of the density, not the velocity gradients.
Let us define the vorticity matrix $\O=(\O_{ij} )$ for the $n$ dimensional vector field $v=(v^1, \cdots, v^n)$ as
$\O_{ij} := \partial_i v^j -\partial_j v^i$, $i,j=1,\cdots, n$.  In the 1-D case we set $\O\equiv 0$. We state our main theorem for the system (\ref{ep1}).
\begin{thm}
Let $n\geq 1$ in the Euler-Poisson equations (\ref{ep1}).   We suppose the initial data $(v_0, \rho_0 )$ satisfies the following condition,
$$
\mathcal{ S}_1=\{a\in \Bbb R^n\,|\, \O_0 (a)=0,\,  -\mathrm{ div}\, v_0 (a)  \geq \sqrt{\frac{2k\rho_0 (a)}{3}} >0\}\neq\emptyset .
$$
Then, the $C^1-$regularity of local classical solution with initial data $v_0, \rho_0$ cannot persist arbitrarily long time.
\end{thm}
\noindent{\bf Remark 1.1 } We observe that in the 1-D case and in the n-D spherically symmetric case the vorticity vanishing condition, $\O _0 (a)=0$ is redundant.\\
\ \\
\noindent{\bf Remark 1.2 } As will be seen on the proof in the next section  the proof does not depend on the specific domain, and the same results hold true for the bounded domain with smooth boundary or the periodic domain.\\
\ \\
Next, we consider the following Euler-Poisson equations with constant background density, $\rb$.
\bb\label{ep2}
 \left\{ \aligned
 &\partial_t \rho+ \mathrm{div}\,(\rho v)=0,\\
 &\partial_t (\rho v)+\mathrm{div}\,(\rho u\otimes u) =-k\rho \nabla \phi,\\
&\Delta \phi =\rho -\rb,\\
& v(x,0)=v_0 (x), \,\, \rho(x,0)=\rho_0 (x).\\
\endaligned\right.
 \ee
The inclusion of constant  density in the Poisson part of the equation has
 more physically realistic meaning in many cases. For example, if we consider the above system in the periodic
unit box $\Bbb T^n=\Bbb R^n/ \Bbb Z^n$, then the natural choice of $\rb$ is
 $$ \rb=\int_{\Bbb T^n} \rho (x,t) dx,
 $$
 which is constant in time thanks to the first equation  of (\ref{ep2}).
\begin{thm}
Let $n\geq 1$ in the Euler-Poisson equations (\ref{ep2}).   We suppose the initial data $(v_0, \rho_0 )$ satisfies the following condition,
$$
\mathcal{S}_2:=\{\mbox{the set of all points $a\in \Bbb R^n$ satisfying (\ref{ep3aa})-(\ref{ep3ab}) below} \}\neq \emptyset
$$
where
\bb\label{ep3aa}
 \O_0 (a)=0, \,\,\, \mathrm{div}\, v_0 (a)<0, \,\,\, \rho_0 (a)\geq \left(\frac12 +\sqrt{\frac{3}{2}} \right)\bar{\rho},
 \ee
 and
 \bb\label{ep3ab}
(\mathrm{ div}\, v_0 (a) )^2 >\max\left\{\frac{k(\rho_0(a) -\frac{\rb}{2}) \left[4 (\rho_0 (a)-\frac{\rb}{2} )^2 -3\rb^2 \right]}{6\rho_0 (a)^2}, \, \frac{3k \rb^4}{4 \rho_0 ^2(a)  (\rho_0(a) -\frac{\rb}{2})}\right\}.
\ee
Then, the $C^1-$regularity of local classical solution with initial data $v_0, \rho_0$ cannot persist arbitrarily long time.
\end{thm}

\section{Proof of  Theorems 1.1 and 1.2}
 \setcounter{equation}{0}

\noindent{\bf Proof of  Theorem 1.1 } We write the first two equations of  (\ref{ep1}) in a more convenient form,
\bb\label{ep4}
 \partial_t \rho+ (v\cdot \nabla)\rho =-\rho \,\mathrm{div}\,v,
 \ee
 and
 \bb\label{ep5}
 \partial_t v+(v\cdot \nabla)v=-k \nabla \phi
 \ee
 respectively.
Taking a partial derivative of (\ref{ep5}),  we obtain the matrix equation,
\bb\label{ep6}
\partial_t V+(v\cdot \nabla )V+V^2 =-k \Phi,
\ee
where we set $V=(\partial_i v^j )$, and $\Phi=(\partial_{i}\partial_j \phi )$. The symmetric and the antisymmetric parts of (\ref{ep6}) are
\bb\label{ep7}
\frac{D}{Dt} \mathcal{D }= -\mathcal{D }^2 -A^2 - k\Phi,
 \ee
and
 \bb\label{ep8}
\frac{D}{Dt} A=-\mathcal{D }A -A\mathcal{D }
 \ee
 respectively, where  we set $\mathcal{D }=\frac12 (V+V^T)$, $A=\frac12 (V-V^T)(=\frac12\O),$  and denoted $\frac{D}{Dt} =\partial_t +(v\cdot \nabla )$.
We consider evolution along the particle trajectory $\{X(a,t)\}$,
 defined by a solution of the following ordinary differential equations,
 $$
 \frac{\partial X(a,t)}{\partial t}=v(X(a,t),t), \, X(a,0)=a\in \mathcal{S}_1,
 $$
where $v(x,t)$ is a classical solution of the system (\ref{ep1}).
From the equation (\ref{ep8}) we have immediately
 $$
 \frac{D}{Dt} |A|\leq 2|\mathcal{D }||A|,
 $$
where we  used the matrix norm,
 $
 |M|:=\sqrt{\sum_{i,j=1}^n M_{ij}^2}.
 $
By Gronwall's lemma we obtain
 \bb\label{ep9}
 |A(X(a,t),t)|\leq |A_0 (a)|\exp\left[ 2\int_0 ^t |\mathcal{D }(X(a,\tau),\tau )|d\tau \right].
 \ee
Since  $A_0 (a)= \frac12\O_0 (a)=0$ for $a\in \mathcal{S}_1$, we have $A(X(a,t),t)=0$ along the particle trajectory $\{ X(a,t)\}$. Hence, taking trace of (\ref{ep7}), we find
$$
\frac{D}{Dt} (Tr (S))=-|\mathcal{D }|^2-k \Delta \phi
$$
along the trajectories $\{ X(a,t)\}$, which taking into account of the fact $\mathrm{div}\, v=Tr(\mathcal{D })$
and the Poisson part of (\ref{ep1}), leads to
\bb\label{ep10}
\frac{D}{Dt} (\mathrm{div}\, v)=-|\mathcal{D }|^2 -k\rho.
 \ee
 From (\ref{ep4}), using (\ref{ep10}), we have
 \bq\label{ep11}
\frac{D^2 \rho}{Dt^2} &=&-\frac{D}{Dt} (\rho \,\mathrm{div}\, v)=-\left(\frac{D \rho }{Dt}\right)\,\mathrm{div}\, v -\rho \frac{D}{Dt} (\mathrm{div}\, v )\n\\
 &=& \rho ( \mathrm{div}\, v )^2+ \rho |\mathcal{D }|^2 + k \rho^2\n \\
 &\geq& k \rho^2.
 \eq
Integrating this inequality along the particle trajectory once, we have
\bb\label{ep12}
\frac{D\rho}{Dt} (X(a,t),t) \geq \left(\frac{D\rho}{Dt} \right)(t=0)  +k \int_0 ^t \rho^2 (X(a,s),s)ds,
\ee
which implies that
 \bb\label{ep12a}
 \frac{D\rho}{Dt}  (X(a,t),t)\geq 0,
 \ee
  if
\bb\label{ep12b}\left(\frac{D\rho}{Dt} \right)(t=0)=-\rho_0 (a)\,\mathrm{div}\, v_0 (a) \geq 0.
\ee
Note that this is the case for $a\in \mathrm{S}_1$.
Hence, multiplying (\ref{ep11}) by $\frac{D \rho}{Dt} $ without changing direction of the inequality, we find
\bb\label{ep13}
\frac12 \frac{D}{Dt} \left(\frac{D\rho}{Dt} \right)^2 \geq \frac{k}{3}\frac{D\rho^3}{Dt} ,
\ee
and integrating this along the particle trajectory again, we have
\bq\label{ep14}
\left(\frac{D\rho}{Dt}  \right)^2 &\geq& \frac{2k}{3} \rho^3  +  \left(\frac{D\rho}{Dt} (t=0)\right)^2 -\frac{2k}{3} \rho_0 ^3(a) \n \\
&\geq & \frac{k}{3} \rho ^3,
\eq
if
\bb\label{ep14a}
\left(\frac{D\rho}{Dt} (t=0)\right)^2 -\frac{2k}{3} \rho_0 ^3(a) =\rho_0 ^2(a)(\mathrm{div}\, v_0 (a))^2 -\frac{2k}{3} \rho_0 ^3(a)\geq 0,
\ee
which  satisfied in our case $a\in \mathcal{S}_1$. The reduced inequality
$$
\frac{D \rho}{Dt} \geq \sqrt{\frac{k}{3}} \rho^{\frac32}
$$
from (\ref{ep14}) can be immediately solved to yield
\bb\label{ep15}
\rho (X(a,t),t)\geq \frac{\rho_0 (a)}{\left( 1-\sqrt{\frac{k \rho_0 (a)}{12}} t\right)^2},
\ee
which shows  that $\rho(X(a,t),t)$ blows-up at a time no latter than $t_* = \sqrt{\frac{12}{k\rho_0 (a)}}$, if the $C^1$ solution  $(v(x,t), \rho(x,t))$ persists until that time.
$\square$\\
\ \\
\noindent{\bf Proof of  Theorem 1.2 } In this case the proof is the same as that of Theorem 1.1 up to the parts before (\ref{ep11}).
Here, instead of (\ref{ep11}) we have
 \bq\label{ep16}
\frac{D^2  \rho }{Dt^2}&=&-\frac{D}{Dt} (\rho \,\mathrm{div}\, v)=-\left(\frac{D  \rho}{Dt} \right)\,\mathrm{div}\, v -\rho \frac{D}{Dt} (\mathrm{div}\, v )\n\\
 &=& \rho ( \mathrm{div}\, v )^2+ \rho |\mathcal{D }|^2 + k \rho (\rho-\bar{\rho} )\n \\
 &\geq& k \rho (\rho-{\bar{\rho}} )=k \left(\rho-\frac{\rb}{2}\right)^2 -\frac{k \rb^2}{4} .
 \eq
 Setting
 $$\theta=\rho-\frac{\rb}{2},
 $$
  we find that
 \bb\label{ep16a}
\frac{D^2 \theta}{Dt^2} =k \theta^2 -\frac{k \rb^2}{4} .
 \ee
Integrating this inequality along the particle trajectory once, we have
\bq\label{ep17}
\lefteqn{\frac{D \theta}{Dt} (X(a,t),t) \geq \left(\frac{D\theta}{Dt} \right)(t=0) -\frac{k\rb^2 t}{4}  +k \int_0 ^t \theta^2 (X(a,s),s)ds}\hspace{.0in}\n \\
 &&=-\rho_0 (a)\, \mathrm{div}\, v_0 (a) -\frac{k\rb^2 t}{4}  +k \int_0 ^t \theta^2 (X(a,s),s)ds.
\eq
This implies that
 \bb\label{ep17a}
 \frac{D \theta }{Dt}(X(a,t),t)> 0,
 \ee
  if
\bb\label{ep17b}\mathrm{div}\, v_0 (a) <0, \,\,\mbox{and}\,\, 0< t \leq - \frac{4\rho_0 (a)\, \mathrm{div}\, v_0 (a)}{k \rb^2}.
\ee
Note that for $a\in \mathcal{S}_2$ the first condition is satisfied, and we will concentrated
on the time interval $[0, t]$ defined  by the second inequality in (\ref{ep17b}).
Hence, multiplying (\ref{ep16a}) by $\frac{D\theta}{Dt} $ without changing direction of the inequality, we find
\bb\label{ep18}
\frac12 \frac{D}{Dt} \left(\frac{D  \theta}{Dt}\right)^2 \geq \frac{k}{3}\frac{D \theta^3}{Dt}  -\frac{k \rb^2 }{4}\frac{D\theta}{Dt} ,
\ee
and integrating along the particle trajectory again, we have
\bq\label{ep19}
\left(\frac{D \theta }{Dt} \right)^2 &\geq& \frac{2k}{3} \theta^3 -\frac{k \rb^2 }{2} \theta +  \left(\frac{D\theta}{Dt} (t=0)\right)^2 -\frac{2k}{3} \theta_0 ^3(a) +\frac{k \rb^2 }{2} \theta_0 (a) \n \\
&=&\frac{2k}{3} \theta^3 -\frac{k \rb^2 }{2} \theta +  \rho_0 ^2(a)(\mathrm{div}\, v_0 (a))^2 -\frac{2k}{3} \theta_0 ^3(a)+\frac{k \rb^2 }{2} \theta_0 (a)\n \\
&\geq&\frac{2k}{3} \theta ^3-\frac{k \rb^2 }{2} \theta =\frac{k}{3} \theta ^3 +\left(\frac{k}{3} \theta ^3-\frac{k \rb^2 }{2} \theta\right)\n \\
&\geq & \frac{k}{3} \theta ^3,
\eq
if
\bb\label{ep19a}
\rho_0 ^2(a)(\mathrm{div}\, v_0 (a))^2 -\frac{2k}{3} \theta_0 ^3(a)+\frac{k \rb^2 }{2} \theta_0 (a)\geq  0,
\ee
and
\bb\label{ep19b}
\frac{k\theta_0 (a)^3}{3}\geq  \frac{k \rb^2 }{2} \theta_0(a).
\ee
Indeed,  from (\ref{ep17a}), we see that (\ref{ep19b}) implies that
$$ \frac{k\theta^3(X(a,t),t)}{3}\geq \frac{k \rb^2 }{2} \theta (X(a,t),t) $$
as long as the classical solution is well-defined.
The conditions (\ref{ep19a}) and (\ref{ep19b}) are satisfied, in turn, if
 \bb\label{ep19c}
 (\mathrm{div}\, v_0 (a))^2 \geq \frac{k(\rho_0(a) -\frac{\rb}{2}) \left[4 (\rho_0 (a)-\frac{\rb}{2} )^2 -3\rb^2 \right]}{6\rho_0 (a)^2},
 \ee
 and
 \bb
 \rho_0 (a) \geq \left(\frac12 +\sqrt{\frac{3}{2}} \right)\bar{\rho}
 \ee
 respectively,
which  follows from the fact $a\in \mathcal{S}_2$. The inequality
$$
\frac{D\theta}{Dt}  \geq \sqrt{\frac{k}{3}} \theta ^{\frac32}
$$
can be solved to yield
\bb\label{ep20}
\theta (X(a,t),t)\geq \frac{\theta_0 (a)}{\left( 1-\sqrt{\frac{k \theta_0 (a)}{12}} t\right)^2},
\ee
which shows  that $\rho(X(a,t),t)=\theta(X(a,t),t) +\frac{\rb}{2}$ blows-up at a time no latter than $t_* = \sqrt{\frac{12}{k\theta_0 (a)}}$, if the $C^1$ solution  $(v(x,t), \rho(x,t))$ persists until that time.
Comparing this with the second inequality of (\ref{ep17b}), in order to get the finite time singularity, it suffices to have
$$
- \frac{4\rho_0 (a)\, \mathrm{div}\, v_0 (a)}{k \rb^2} > \sqrt{\frac{12}{k\theta_0 (a)}}=t_*,
$$
which is true, if
\bb\label{ep20a}
(\mathrm{div}\, v_0 (a) )^2 > \frac{3k \rb^4}{4 \rho_0 ^2(a)  (\rho_0(a) -\frac{\rb}{2})}.
\ee
Both of the conditions (\ref{ep19c}) and (\ref{ep20a}) are satisfied for $a\in \mathcal{S}_2$.
$\square$\\
\section{Remarks on the Polytropic Euler equations}
 \setcounter{equation}{0}

 In this section we show that a similar argument as in the above section can be used to show a finite time blow-up result for the following Euler equations for polytropic gas flows in a domain in $\Bbb R^n$.

 \bb\label{21}
 \left\{ \aligned
 &\partial_t \rho+ \mathrm{div}\,(\rho v)=0,\\
 &\partial_t (\rho v)+\mathrm{div}\,(\rho v\otimes v) =-\nabla p,\\
&\partial_t S+(v\cdot \nabla )S =0,\\
& p=p(\rho, S)=\rho^\g e^S, \quad (\g >1)\\
& v(x,0)=v_0 (x), \,\, \rho(x,0)=\rho_0 (x),\,\, S(x,0)=S_0.\\
\endaligned\right.
 \ee
 In the above $p=p(\rho, S)$ is the pressure, and $S$ is the entrophy.
 For both of the physical and mathematical backgrounds of the system (\ref{21}) including the local well-posedness of classical solutions we refer the monograph \cite{maj} and the review articles \cite{ che,chr}.
For the blow-up problem of the system (\ref{21}) for $n=3$, in particular, there is already  well-known result due to Sideris in \cite{sid}(see also \cite{che, chr} and the references therein for the other related results on this problem). In this section we just would like to show another aspect of the blow-up problem of the Euler system, represented in Theorem 3.1 below.\\
 As in the previous section we rewrite the system (\ref{21}) in a more convenient form,
  \bb\label{22}
 \left\{ \aligned
 &\partial_t \rho+ (v\cdot \nabla )\rho= -\rho \,\mathrm{div} \, v,\\
 &\partial_t  v+(v\cdot \nabla) v  =-\nabla p,\quad p=e^S \rho ^\g\\
&\partial_t S+ (v\cdot \nabla ) S =0,\\
& v(x,0)=v_0 (x), \,\, \rho(x,0)=\rho_0 (x),\,\, S(x,0)=S_0.\\
\endaligned\right.
 \ee
We use the matrix notations, $\mathcal{D}_{ij}= \frac12 (\partial_i v^j +\partial_j v^i)$,
$A_{ij}= \frac12 (\partial_i v^j -\partial_j v^i)$, and $P_{ij}=\partial_i\partial_j p$ as previously.
We also denote $\{ (\lambda_k (x,t),e_k (x,t))\}$ for the pairs of eigenvalue and the normalized eigenvector of the deformation tensor $\mathcal{D}$.
 In this section we shall prove the following theorem.
\begin{thm}
Let $\g\geq 2$ be given. We suppose the initial data $(v_0, \rho_0 , S_0 )$  for (\ref{22}) satisfies the following condition,
$$
\mathcal{S}_3:=\{\mbox{the set of all points $a\in \Bbb R^n$ satisfying (\ref{22a})-(\ref{22b}) below}\}\neq \emptyset
$$
where
\bb\label{22a}
 \O_0 (a)=0,
 \ee
 \bb\label{22ab}
  D^k S_0 (a)=0\, \forall  \,k\in\{1,2\}\,\mbox{and} \,\, D^m \rho_0 (a)=0
\, \forall \, m\in \{0,1,2\},
 \ee
 and
 \bb\label{22b}
\exists j\in \{1,\cdots, n\}\,\,\mbox{such that}\,\,\lambda_j (a,0)<0
\ee
Then, the  $C^3-$regularity of solution to the system (\ref{22}) cannot persist arbitrary long time.
\end{thm}
 {\bf Proof } Taking the first and the second partial derivatives of the first equation of (\ref{22}), we find that
 \bb\label{23}
 \frac{D} {Dt}\partial_j \rho=-(\partial_j v\cdot \nabla ) \rho -\partial_j \rho \,\mathrm{div}\, v -\rho  \mathrm{div}\, \partial_j v,
 \ee
 and
 \bq\label{24}
\frac{D} {Dt}\partial_j\partial_k \rho&=&-(\partial_j\partial_k  v\cdot \nabla ) \rho -(\partial_j v\cdot \nabla ) \partial_k \rho-(\partial_k v\cdot \nabla )\partial_j \rho-
 \partial_j \partial_k \rho \,\mathrm{div}\, v \n \\
 &&\quad -\partial_j \rho \,\mathrm{div}\, \partial_k v -\rho \, \mathrm{div}\, \partial_j \partial_k v-\partial_k\rho  \,\mathrm{div}\, \partial_j v
 \eq
 respectively.
Combining (\ref{23}) and (\ref{24}) with the first equation of (\ref{22}), we find that
 $$
 \frac{D} {Dt}(\rho +|D\rho |+|D^2\rho|)\leq 2(|Dv|+|D^2v |+|D^3 v|) (\rho +|D\rho |+|D^2\rho|),
$$
 which, after integration along the particle trajectories, provides us with
 \bq\label{25}
 \lefteqn{(\rho +|D\rho |+|D^2\rho|)(X(a,t),t)\leq }\hspace{.0in}\n \\
  &&\leq (\rho_0 +|D\rho_0 |+|D^2\rho_0|)(a)\exp\left[2\int_0 ^t
  (|Dv|+|D^2v |+|D^3 v|) (X(a,s),s)ds\right].\n \\
 \eq
 Similarly, from the third equation of (\ref{22}) we obtain the estimate,
  \bq\label{26}
 \lefteqn{(|D S |+|D^2S|)(X(a,t),t)\leq }\hspace{.0in}\n \\
  &&\leq (|D S_0 |+|D^2 S_0|)(a)\exp\left[2\int_0 ^t
  (|Dv|+|D^2v |) (X(a,s),s)ds\right].\n \\
 \eq
 The symmetric and the antisymmetric parts of the velocity gradient matrix equations, obtained from the second equation of (\ref{22}) are
 \bb\label{27}
 \frac{D} {Dt}\mathcal{D}=-\mathcal{D}^2 -A^2 -P,
 \ee
 and
 \bb\label{28}
 \frac{D} {Dt} A=-\mathcal{D}A -A \mathcal{D}
 \ee
respectively.  In particular, from (\ref{28}) we obtain easily that
\bb\label{29}
|A (X(a,t),t)|\leq |A_0(a)|\exp\left[2\int_0 ^t
  |Dv (X(a,s),s)|ds\right].
\ee
We observe the pointwise estimate for the Hessian of the pressure,
\bb\label{30}
|P|\leq C e^S (|DS|+|D^2 S|)(\rho +|D\rho|+|D^2 \rho |),
\ee
where $C=C(\rho )$, and $\lim\sup_{\rho \to 0}|C( \rho )| <\infty $ thanks to the assumption $\g \geq 2$.
The estimates (\ref{25}), (\ref{26}), (\ref{29}) and (\ref{30}) shows that
\bb\label{31}
A(X(a,t),t)=P (X(a,t),t)=0
\ee
if $a\in \mathcal{S}_3$, which we suppose from now on.
Hence, along the particle trajectories $\{ X(a,t)\}$, the equation (\ref{27}) reduces to
\bb\label{32}
 \frac{D} {Dt}\mathcal{ D}=- \mathcal{D}^2.
 \ee
 Let $(\lambda_j, e_j) $ be an eigenvalue-normalized eigenvector pair of $\mathcal{D}$ such that $\lambda_j(a,0)<0$.
Then, taking operation $\frac{D}{Dt}$ on $\mathcal{D}e_j=\lambda_j e_j$, we have
 $$ \left(\frac{D }{Dt}\mathcal{D}\right)e_j + \mathcal{D} \frac{D e_j }{Dt}=\frac{D \lambda_j }{Dt} e_j +\lambda_j \frac{D e_j }{Dt}.$$
 Using (\ref{32}), this reduces to
  \bb\label{32a}
  -\lambda_j^2 e_j + \mathcal{D} \frac{D e_j }{Dt}=\frac{D \lambda_j }{Dt} e_j +\lambda_j \frac{D e_j }{Dt}.
  \ee
  Multiplying $e_j$ on the both sides of (\ref{32a}), we find
  \bb
  \frac{D \lambda_j }{Dt}=-\lambda_j^2,
  \ee
which can be solved along the particle trajectories $\{ X(a,t)\}$ to yield
  $$
  \lambda_j  (X(a,t),t)= \frac{  \lambda _{j} (a,0)}{1+\lambda _{j} (a,0) t},
  $$
  and shows that
  $$ \lambda_j (X(a,t),t)\to -\infty \quad \mbox{as}\quad t\to t_*=\frac{-1}{\lambda_{j}(a,0)}.
   $$
    $\square$
\[ \mbox{\bf Acknowledgements} \]
This research was done while the author was visiting University of Nice. The author would like to thank to Y. Brenier  for his hospitality and for stimulating discussions with useful suggestions. He would like to thank also to E. Tadmor for informing him of the recent results on the subject with references. This research was supported partially by KRF Grant(MOEHRD, Basic Research Promotion Fund).

  \end{document}